\documentclass[9pt,letterpaper,dvips]{article}
 \usepackage{amssymb,latexsym,amsmath,amsthm}

 \renewcommand{\div}{\mathop{\mathrm{div}}\nolimits}

\newtheorem{thm*}{Theorem}

\newtheorem{thm}{Theorem}

\newtheorem{remark}{Remark}

\newtheorem{dfn}{Definition}
\newtheorem{prop}{Proposition}

\newtheorem{conj}{Conjecture}

\newtheorem{lemma}{Lemma}

\newtheorem{rem}{Remark}

\begin{document}
\date{}

\title{De Giorgi type results for elliptic systems}
\author{Mostafa Fazly\thanks{Research partially supported by a University Graduate Fellowship at the University of British Columbia, under the supervision of the second-named author.} \quad and \quad Nassif  Ghoussoub\thanks{Partially supported by a grant
from the Natural Sciences and Engineering Research Council of Canada.  }
\quad   
\\
\small Department of Mathematics,\\
\small University of British Columbia, \\
\small Vancouver BC Canada V6T 1Z2 \\\\
\small {\tt fazly@math.ubc.ca}\\
\small {\tt nassif@math.ubc.ca} \\\\
\date{Revised April 15, 2012}\\
}
\maketitle

\vspace{3mm}

\vspace{3mm}

\begin{abstract} We consider the following elliptic system 
\begin{eqnarray*}
 \Delta u =\nabla H (u) \ \ \text{in}\ \ \mathbf{R}^N,
  \end{eqnarray*}
where $u:\mathbf{R}^N\to \mathbf{R}^m$ and $H\in C^2(\mathbf{R}^m)$, and prove, under various conditions on the nonlinearity $H$ that, at least in low dimensions, a solution $u=(u_i)_{i=1}^m$ is necessarily one-dimensional whenever each one of its components $u_i$ is monotone in one direction. Just like in the proofs of the classical De Giorgi's conjecture in dimension $2$ (Ghoussoub-Gui) and in dimension $3$ (Ambrosio-Cabr\'{e}), the 
key step is a Liouville theorem for linear systems. We also give an extension of a geometric Poincar\'{e} inequality to systems and use it to establish De Giorgi type results for stable solutions as well as  additional rigidity properties stating that the gradients of the various components of the solutions must be parallel. We introduce and exploit  the concept of {\it an orientable system}, which seems to be key for dealing with systems of three or more equations. For such systems, the notion of a stable solution in a variational sense coincide with the pointwise (or spectral) concept of stability.

\end{abstract}


\section{Introduction}

In 1978, Ennio De Giorgi proposed the following conjecture.

\begin{conj}\label{conj1} 
Suppose that $u$ is an entire solution of the Allen-Cahn equation \begin{equation}\label{allen}
\Delta u+u-u^3 =0	 \quad \text{on}\ \  \mathbf{R}^N
\end{equation}
satisfying $|u(\mathbf{x})| \le 1$, $\frac{\partial u}{\partial x_N} (\mathbf{x}) > 0$ for $\mathbf{x} = (\mathbf{x}',x_N) \in \mathbf{R}^N$.	
Then, at least in dimensions $N\le 8$ the level sets of $u$ must be hyperplanes, i.e. there exists $g \in C^2(\mathbf{R})$ such that $u(\mathbf{x}) = g(a\mathbf{x}' - x_N)$, for some fixed $a \in \mathbf{R}^{N-1}$.
\end{conj}
The first positive result on the De Giorgi conjecture was established in 
1997 by Ghoussoub and Gui \cite{gg1} 
for dimension $N= 2$. 
Their proof used the following linear Liouville type theorem for elliptic equations in divergence form, which (only) holds in dimensions $1$ and $2$ (\cite{bbg, gg1}). If $\phi>0$, then any solution $\sigma$ of 
\begin{equation}
\div(\phi^2 \nabla \sigma)=0, 
\end{equation}
such that $\phi\sigma$ is bounded, is necessarily constant. This result is then applied to 
 the ratio $\sigma := \frac{   {\frac{\partial u}{\partial x_1} } } {\frac{\partial u}{\partial x_2}  }$ to conclude in dimension 2. 
Ambrosio and Cabr\'{e} \cite{ac} extended the result to dimension $N= 3$ by noting  that for the linear Liouville theorem to hold,  it suffices that 
\begin{equation}
\int_{B_R}\phi^2\sigma^2 \leq CR^2,
\end{equation} 
and  by proving that any solution $u$ satisfying $\partial_N u>0$ satisfies the energy estimate 
\begin{equation}
\int_{B_R}\phi^2\sigma^2 \leq CR^{N-1}.
\end{equation} 
The conjecture remains open in dimensions $4\leq N \leq 8$. However, 
Ghoussoub and Gui also showed in \cite{gg2} that it is true for $N = 4$ or $N = 5$ for solutions that satisfy certain antisymmetry conditions, and Savin \cite{s} established its validity   for $4 \le N \le 8$ under the following additional natural hypothesis on the solution,
 \begin{equation}\label{asymp}
\lim_{x_N\to\pm\infty } u(\mathbf{x}',x_N)\to \pm 1.
\end{equation}
Unlike the above proofs in dimensions $N\leq 5$, the proof of Savin is nonvariational and does not use a Liouville type theorem. Our proofs below for analogous results corresponding to systems are more in the spirit of Ghoussoub-Gui and Ambrosio-Cabr\'{e},  which mostly rely on notions of stability and on an interesting linear Liouville theorem that is suitable for non-linear elliptic systems of the following type: 
\begin{eqnarray}
\label{main}
 \Delta u =\nabla H (u) \ \ \text{in}\ \ \mathbf{R}^N,
  \end{eqnarray}
where $u:\mathbf{R}^N\to \mathbf{R}^m$, $H\in C^2(\mathbf{R}^m)$ and $\nabla H (u)=(H_{u_i}(u_1, u_2,...u_m))_{i}$. The notation $H_{u_i}$ is for the partial derivative $\frac{\partial H}{\partial u_i}$.

 \begin{dfn} \label{weak} \rm 
 We shall say that the system (\ref{main}) (or the non-linearity $H$) is {\it orientable}, if there exist nonzero functions $\theta_k\in C^1(\mathbf{R}^N)$, $k=1,\cdots,m$, which do not change sign, such that for all $i,j$ with $1\leq  i<j\leq m$, we have 
 \begin{equation}
 \hbox{$ H_{u_iu_j}\theta_i(x)\theta_j(x)\le 0$ \, for all $x\in\mathbf{R}^N$.}
  \end{equation} 
\end{dfn}
Note that the above condition on the system means that none of the mixed derivative $H_{u_iu_j}$  changes sign. It is clear that a system consisting of two equations  (i.e., $m=2$)  is always orientable as long as $H_{u_1u_2}$ does not change sign.  On the other hand, if $m=3$, then the system  (\ref{main}) cannot be orientable if, for example, all three mixed derivatives $H_{u_iu_j}$ with $i<j$ are positive. This concept of "orientable system" seems to be the right framework for dealing with systems of three or more equations. We shall see for example, that for such systems, the notions of variational stability and pointwise stability coincide.

We shall consider solutions of (\ref{main}) whose components $(u_1, u_2,...u_m)$ are strictly monotone in the last variable $x_N$. However, and in contrast to the case of a single equation, the various components  need not be all increasing (or decreasing).  This leads us to the following definition of monotonicity.
  
 \begin{dfn} \rm
Say that a solution $u=(u_k)_{k=1}^m$ of (\ref{main}) is {\it $H$-monotone} if the following hold:
\begin{enumerate}
 \item For every $i\in \{1,..., m\}$, $u_i$ is strictly monotone in the $x_N$-variable (i.e., $\partial_N u_i\neq 0$).

\item  For $i<j$, 
we have 
  \begin{equation}
\hbox{$H_{u_iu_j} \partial_N u_i(x) \partial_N u_j (x)\le 0$  for all $x\in\mathbf{R}^N$.}
\end{equation}

\end{enumerate}
We shall then write $ I\cup J=\{1,...,m\}$, where 
 \begin{eqnarray}\label{mono}
\hbox{$\partial_N u_i>0 > \partial_N u_j$ for $i\in I$ and $j\in J$. }
\end{eqnarray}

\end{dfn}
It is clear that the mere existence of an $H$-monotone solution for (\ref{main}) implies that the system  is orientable, as it suffices to use  $\eta_i=\partial_N u_i$. We now recall two notions of stability that will be considered in the sequel.

\begin{dfn} \rm 
A solution $u$ of the system (\ref{main}) on a domain $\Omega$ is said to be
\begin{enumerate}
\item [(i)]  {\it stable}, if the second variation of the corresponding energy functional is nonnegative, i.e., if 
\begin{equation}\label{stability} 
\sum_{i} \int_{\Omega}  |\nabla \zeta_i|^2  +\sum_{i,j} \int_{\Omega} H_{u_iu_j} \zeta_i\zeta_j \ge 0,
 \end{equation}
for every $\zeta_k \in C_c^1(\Omega), k=1,...,m$.

\item [(ii)] {\it pointwise stable}, if there exist $(\phi_i)_{i=1}^m$ in $C^1(\Omega)$ that do not change sign and $\lambda\ge 0$ such that 
\begin{equation}\label{pointstability} 
\hbox{$\Delta \phi_i=\sum_{j} H_{u_i,u_j} \phi_j - \lambda \phi_i$ in $\Omega$\,  for all $i=1,...,m$,}
 \end{equation}
and $H_{u_i,u_j} \phi_j \phi_i\le 0$ for $1\leq i<j\leq m$.
\end{enumerate}
 \end{dfn}
 Note that a system that possesses a pointwise stable solution is necessarily orientable. On the other hand, we shall prove in Section \ref{sec1} that for solutions of orientable systems, the two notions of stability are equivalent. The main focus of this paper is to provide some answers to the following conjecture, which extends the one by De Giorgi on Allen-Cahn equations to more general systems.

\begin{conj}\label{conj2} 
Suppose $u=(u_i)_{i=1}^m$ is an $H$-monotone bounded entire solutions of the system (\ref{main}), then at least in dimensions $N\le 8$, the level sets of each component $u_i$ must be a hyperplane. 
\end{conj}

We are extremely grateful to an anonymous referee for several important comments, including the suggestion that our methods could also lead to a solution of the conjecture in dimension $3$, without the additional assumptions that we had used in the first version of this paper.

\section{A linear Liouville theorem for systems and first applications}\label{sec1}

The following Liouville theorem plays a key role in this paper. Note that for the case $m=1$, this type of Liouville theorem was noted by Berestycki, Caffarelli and Nirenberg in \cite{bcn} and used by Ghoussoub-Gui \cite{gg1} and later by Ambrosio and Cabr\'{e} \cite{ac}  to prove the De Giorgi conjecture in dimensions two and three. Also, Ghoussoub and Gui in \cite{gg2} used a slightly stronger version to show that the De Giorgi's conjecture is true in dimensions four and five for a special class of solutions that satisfy an antisymmetry condition.

\begin{prop}\label{liouville} Assume that $\phi_i \in L^{\infty}_{loc}(\mathbf{R}^N) $ are such that $\phi^2_i >0 $ a.e.,  and $\sigma_i \in H^1_{loc}(\mathbf{R}^N)$ satisfy 
\begin{equation}\label{liouassum}
 \sum_{i=1}^{m}\int_{B_{2R}\setminus B_R}  \phi_i^2\sigma_i^2\le C R^2.
 \end{equation}
 If $(\sigma_i)_{i=1}^m$ are solutions of 
\begin{eqnarray}
\label{div}
\div(\phi_i^2 \nabla \sigma_i) +\sum_{j=1}^{m} h_{i,j}(x) f(\sigma_i-\sigma_j)+k_i(x)g(\sigma_i)= 0   \ \ \text{in}\ \ \mathbf{R}^N\ \ \ \ \text{for} \ \ i=1,...,m ,
  \end{eqnarray}
where $0\ge h_{i,j},k_i\in L_{loc}^1(\mathbf{R}^N)$, $h_{i,j}=h_{j,i}$ and $f,g\in L_{loc}^1(\mathbf{R})$ are odd functions such that $f(t),g(t)\ge 0$ for $t\in\mathbf R^+$.  Then, for all  $i=1,...,m$, the functions $\sigma_i$ are constant.
\end{prop}
\textbf{Proof:}  Multiply both sides of (\ref{div}) by $\sigma_i\zeta_R^2$  where $\zeta_R\in C^1_c(\mathbf{R}^N)$ with $0\le\zeta_R\le1$ being the following test function;
 $$\zeta_R(x)=\left\{
                      \begin{array}{ll}
                        1, & \hbox{if $|x|<R$,} \\
                        0, & \hbox{if $|x|>2R$,} 
                                                                       \end{array}
                    \right.$$
where $||\nabla\zeta_R||_{\infty}\le R^{-1}$. By integrating by parts, we get 
\begin{eqnarray*}
\label{}
\int_{B_{2R}}  \phi_i^2 |\nabla \sigma_i|^2\zeta_R^2 + 2\int_{B_{2R}}  \phi^2\nabla\sigma_i\cdot\nabla \zeta_R \zeta_R \sigma_i- \int_{B_{2R}} \sum_{j} h_{i,j}(x) f(\sigma_i-\sigma_j)\sigma_i \zeta_R^2- \int_{B_{2R}}k_i(x)g(\sigma_i)\sigma_i \zeta_R^2&=& 0.
  \end{eqnarray*}
Summing the above identity over $i$, we get 
\begin{equation*}
\sum_{i=1}^{m} \int_{B_{2R}} \phi_i^2 |\nabla \sigma_i|^2\zeta_R^2 = - 2 \sum_{i=1}^{m} \int_{B_{2R}\setminus B_R}  \phi_i^2\nabla\sigma_i\cdot\nabla \zeta_R \zeta_R \sigma_i  + \int_{B_{2R}}  I(x) \zeta_R^2+\int_{B_{2R}}  J(x) \zeta_R^2,
\end{equation*}
where $$I(x):= \sum_{i,j} h_{ij}(x) \sigma_i f(\sigma_i-\sigma_j) \ \ \text{and}\ \ J(x):= \sum_{i} k_{i}(x)g(\sigma_i)\sigma_i .$$ 
Note that 
\begin{eqnarray*}
I(x)&=& \sum_{i,j} h_{ij}(x) \sigma_i f(\sigma_i-\sigma_j)\\&=& \sum_{i< j}  h_{i,j} \sigma_i f(\sigma_i-\sigma_j) + \sum_{i> j} h_{i,j}  \sigma_i  f(\sigma_i-\sigma_j)  \\&=&\sum_{i< j} h_{i,j}  \sigma_i f(\sigma_i-\sigma_j)  + \sum_{i< j}h_{i,j} \sigma_j f(\sigma_j-\sigma_i)   \ \ \text{since} \ \ h_{ij}=h_{ji}\\&=&\sum_{i< j} h_{i,j} (\sigma_i-\sigma_j)f(\sigma_i-\sigma_j)  \ \ \text{because $f$ is odd}.
  \end{eqnarray*}
Since $h_{i,j}(x)\le 0$ and $(\sigma_i-\sigma_j)f(\sigma_i-\sigma_j)\ge0$ for all $i,j$, we have $I(x)\le 0$. Similarly, $J(x)\le 0$. Therefore, for $0<\epsilon<1$, 
 we get Young's inequality that  
\begin{eqnarray}\label{decay}
\sum_{i=1}^{m}  \int_{B_{2R}}  \phi_i^2 |\nabla \sigma_i|^2\zeta_R^2&\le& 2 \sum_{i=1}^{m}  \int_{B_{2R}\setminus B_R} \phi_i^2|\nabla\sigma_i| |\nabla \zeta_R|  \zeta_R \sigma_i \nonumber\\ & \le & \epsilon \sum_{i=1}^{m}  \int_{B_{2R}\setminus B_R} \phi_i^2 |\nabla \sigma_i|^2\zeta_R^2 +C_\epsilon\sum_{i=1}^{m}   \int_{B_{2R}\setminus B_R} (\sigma_i\phi_i)^2 |\nabla \zeta_R|^2. \end{eqnarray}
By assumption (\ref{liouassum}) we see that 
\[
\hbox{$\sum_{i=1}^{m}  \int_{\mathbf{R}^N}  \phi_i^2 |\nabla \sigma_i|^2\zeta_R^2  <\infty$\quad 
and \quad 
$\sum_{i=1}^{m}  \int_{\mathbf{R}^N} (\sigma_i\phi_i)^2 |\nabla \zeta_R|^2 <\infty. $}
\]
Estimate (\ref{decay}) then yields
\begin{eqnarray*}
\sum_{i=1}^{m}  \int_{\mathbf{R}^N}  \phi_i^2 |\nabla \sigma_i|^2\zeta_R^2 =0, 
\end{eqnarray*}
which means that $\sigma_i$ for all $i=1,...,m$ must be constant.\hfill $ \Box$\\

Our first application is the following extension of a recent result by Berestycki, Lin,  Wei and  Zhao \cite{blwz} who considered a system of $m=2$ equations and the nonlinearity $H(t,s)=\frac{1}{2}t^2s^2$, which also appear as a limiting elliptic system arising in phase separation for
multiple states Bose-Einstein condensates.

\begin{thm} Suppose the nonlinearity $H$ satisfies the condition:
\begin{equation}
u_iH_{u_i} \ge 0 \ \ \ \text{for all $1\le i\le m$}. 
\end{equation}
Then, any pointwise stable solution $u$ of the system (\ref{main}), which satisfies
\begin{equation}
\sum_{i}\int_{B_{2R}\setminus B_R} u_i^2\le C R^4,
\end{equation}
is necessarily one-dimensional.
 \end{thm}
\textbf{Proof:} Note that we do not assume here that $u$ is bounded solution.  Multiply both sides of (\ref{main}) with $u_i\zeta^2$ to get 
\begin{eqnarray}
\label{supersol}
\hfill u_i \Delta u_i \zeta^2= u_i H_{u_i} \zeta^2&\ge& 0  \ \ \text{in}\ \ \mathbf{R}^N. 
  \end{eqnarray}
An integration by parts yields,  
\begin{eqnarray}
\label{}
\int_{B_R} |\nabla  u_i|^2\zeta^2&\le& 2 \int_{B_R} |\nabla u_i||\nabla \zeta| u_i\zeta. 
  \end{eqnarray}
Now, use the same test function as in the proof of Proposition \ref{liouville} to obtain
\begin{equation}\label{decayn-2}
\sum_{i}\int_{B_R} |\nabla  u_i|^2\le C R^{-2} \sum_{i} \int_{B_{2R}\setminus B_R}  u_i^2 \le C R^2.
  \end{equation}
Since $u$ is a pointwise stable solution of (\ref{main}), there exist eigenfunctions $(\phi_i)_i$ such that 
\begin{equation}\label{} 
\Delta \phi_i= \sum_{j} H_{u_iu_j}\phi_j -\lambda \phi_i \ \ \text{in}\ \ \mathbf{R}^N,
  \end{equation}
where $\phi_i$ does not change sign, $H_{u_iu_j}\phi_i(x)\phi_j(x)\le 0$ and $\lambda\ge 0$. For any fixed $\eta=(\eta',0)\in \mathbf{R}^{N-1}\times\{0\}$, define $\psi_i:=\nabla u_i\cdot\eta$  and observe that  $\psi_i$ satisfies the following equation
\begin{eqnarray}
 \hfill \Delta \psi_i=  \sum_{j} H_{u_iu_j}\psi_j  \ \ \text{in}\ \ \mathbf{R}^N.
  \end{eqnarray}
It is straightforward to see that $\sigma_i:=\frac{\psi_i}{\phi_i}$  is a solution of system (\ref{div}) with  $h_{i,j}(x)=H_{u_iu_j}\phi_i(x)\phi_j(x)$, $k_i(x)=-\lambda \phi_i^2$ and $f,g$ equal the identity. 
Apply now Proposition \ref{liouville} to deduce that $\sigma_i$ is constant for every $i=1,...,m$, which clearly yields our claim. 
\hfill $ \Box$

\begin{remark}
Note that in the case where $m=2$ and $H(t,s)=\frac{1}{2}t^2s^2$, the above theorem yields that any positive solution $(u,v)$ of the corresponding system (\ref{main}), which satisfies the growth assumption $u(x),v(x)=O(|x|^k)$ and such that $\partial_N u>0$, $\partial_N v<0$, 
is  necessarily one-dimensional provided $N\le 4-2k$. Note that Noris et al. \cite{BTT} have recently shown that a solution such that $ u(x),v(x)\leq C(1+|x|^\alpha)$ is necessarily constant.
\end{remark}

We can also deduce the following Liouville theorem for bounded solutions of (\ref{main}) with general non-positive nonlinearities. The approach to this Liouville theorem seems to be new, even for single equations. It is worth comparing to the general results of Nedev \cite{Ned} and Cabr\'{e} \cite{c} regarding the regularity of stable solutions of semilinear equations with general nonlinearities up to dimension four.

\begin{thm} Suppose $H$ is a nonlinearity verifying 
\begin{equation}
\hbox{$H_{u_i}\le 0$ for all $i=1,..., m$.}
\end{equation} 
If the dimension $N\le 4$, then any bounded pointwise stable solution of the system (\ref{main}) is necessarily constant.   
 \end{thm}
\textbf{Proof:} Multiply both sides of system (\ref{main}) with $(u_i-||u_i||_{\infty})\zeta^2$. Since $ H_{u_i} (u_i-||u_i||_{\infty})\ge 0$ we have 
\begin{eqnarray}
\label{supersol}
\hfill \Delta u_i (u_i-||u_i ||_\infty)\zeta^2&\ge& 0  \ \ \text{in}\ \ \mathbf{R}^N.
  \end{eqnarray}
After an integration by parts, we end up with 
\begin{eqnarray}
\label{}
\int_{B_R} |\nabla  u_i|^2\zeta^2&\le& 2 \int_{B_R} |\nabla u_i||\nabla \zeta |(||u_i ||_\infty-u_i)\zeta \ \ \text{for all $1\le i\le m$.}
  \end{eqnarray}
Using Young's inequality and adding we get 
\begin{equation}\label{decayn-2}
\sum_{i} \int_{B_R} |\nabla  u_i|^2 \le R^{N-2}.
  \end{equation}
As in the preceding theorem, one can apply Proposition \ref{liouville} to quotients of partial derivatives to obtain that each $u_i$ is one dimensional solutions as long as $N\le 4$. Note now that $u_i$ is a bounded solution for (\ref{supersol}) in dimension one, and the corresponding decay estimate (\ref{decayn-2}) now implies that  $u_i$ must be constant for all $1\le i\le m$.  \hfill $ \Box$\\

We now show that stability and pointwise stability are equivalent for solutions of orientable elliptic systems. 

\begin{lemma} A $C^2$-function is a pointwise stable solution of the system (\ref{main}) if and only if 
 it is a stable solution and the system is orientable.
\end{lemma}

\textbf{Proof:} Assume first that $u$ is a pointwise stable solution for (\ref{main}). It is clear that the system is then obviously orientable. In order to show that $u$ is a stable solution, we consider  test functions  $\zeta_i \in C_c^1(\mathbf{R}^N)$ and multiply both sides of (\ref{pointstability}) with $\frac{\zeta_i^2}{\phi_i}$ to obtain
\begin{eqnarray*}
&&-\int_{\mathbf{R}^N}  |\nabla\phi_i|^2 \frac{\zeta_i^2}{\phi^2_i}+ 2 \int_{\mathbf{R}^N} \nabla \phi_i\cdot\nabla \zeta_i\ \frac{\zeta_i}{\phi_i}+\sum_{j}\int_{\mathbf{R}^N} H_{u_iu_j} \frac{\phi_j}{\phi_i}\zeta_i^2-\lambda \zeta_i^2=0.
  \end{eqnarray*}
 By applying Young's inequality for the first two terms and taking sums, we get 
 \begin{eqnarray*}
\sum_{i} \int_{\mathbf{R}^N}  |\nabla\zeta_i|^2 +\sum_{i,j} \int_{\mathbf{R}^N} H_{u_iu_j}\frac{\phi_j}{\phi_i}\zeta_i^2\ge\lambda \sum_{i}  \int_{\mathbf{R}^N} \zeta_i^2\ge 0.
  \end{eqnarray*} 
  Note now that 
     \begin{eqnarray*}
  \sum_{i,j} H_{u_iu_j}\frac{\phi_j}{\phi_i}\zeta_i^2 &=& \sum_{i} H_{u_iu_i} \zeta_i^2 + \sum_{i\neq j} H_{u_iu_j}\frac{\phi_j}{\phi_i}\zeta_i^2 \\&=& \sum_{i} H_{u_iu_i} \zeta_i^2 + \sum_{i < j} H_{u_iu_j}\frac{\phi_j}{\phi_i}\zeta_i^2 + \sum_{i> j} H_{u_iu_j}\frac{\phi_j}{\phi_i}\zeta_i^2 \\&= & \sum_{i} H_{u_iu_i} \zeta_i^2 + \sum_{i < j} H_{u_iu_j}\frac{\phi_j}{\phi_i}\zeta_i^2 + \sum_{i< j} H_{u_iu_j}\frac{\phi_i}{\phi_j}\zeta_j^2 \\ &=& \sum_{i} H_{u_iu_i} \zeta_i^2 + \sum_{i < j} H_{u_iu_j}  (\phi_i \phi_j)^{-1} \left( \phi_j^2\zeta_i^2 +\phi_i^2\zeta_j^2\right) \\&\le &  \sum_{i} H_{u_iu_i} \zeta_i^2 + 2\sum_{i < j} H_{u_iu_j} \zeta_i\zeta_j \ \ \text{since $H_{u_iu_j}  (\phi_i \phi_j)^{-1} \le 0$} \\&=&\sum_{i,j} H_{u_iu_j} \zeta_i\zeta_j, 
  \end{eqnarray*} 
  which finishes the proof.
  
   For the reverse implication, we assume the system orientable and consider a stable solution $u$. We shall follow ideas  of Ghoussoub-Gui in \cite{gg1} (see also Berestycki-Caffarelli-Nirenberg in \cite{bcn}) to show that $u$ is pointwise stable. 
   
Define for each $R>0$, 
  \begin{equation}
  \lambda_1(R) := \min_{(\zeta_i)_{i=1}^{m}\in H_0^1(B_R(0)) \setminus\{0\}}  \left\{ \sum_{i} \int_{B_R(0)}  |\nabla \zeta_i|^2 +\sum_{i,j} \int_{B_R(0)} H_{u_i,u_j} \zeta_i\zeta_j , \ \ \ \sum_{i} \int_{B_R(0)} \zeta_i^2=1  \right\}.
  \end{equation}
Since $u$ is a stable solution, we have that $\lambda_1(R)\ge 0$ and there exist eigenfunctions $\zeta_i^R$ such that 
  \begin{equation}\label{zetar}
  \left\{
                      \begin{array}{ll}
                                             \Delta \zeta_i^R= \sum_{j} H_{u_iu_j} \zeta^R_j- \lambda_1(R) \zeta_i^R, & \hbox{if $|x|< R$,} \\
                       \zeta_i^R=0, & \hbox{if $|x|= R$.}
                                                                       \end{array}
                    \right.\end{equation}
Since the system  is orientable, there exists $(\theta_k)_{k=1}^{m}$ such that $ H_{u_iu_j}\theta_i\theta_j\le 0.$ We can then use the signs of the $\theta_k$'s to assign signs for the eigenfunctions $(\zeta_k^R)_k$ so that they satisfy 
\begin{equation}\label{zetamono}
 H_{u_iu_j}\zeta_i^R\zeta_j^R\le 0.
 \end{equation} 
For that it suffices to replace $\zeta_i^R$ with $sgn(\theta_i)|\zeta_i^R|$ if need be. We can also normalize them so that 
  \begin{equation}\label{normal}
  \sum_{k}  |\zeta_k^R(0)|=1.
  \end{equation}
  Note that $\lambda_1(R) \downarrow \lambda_1 \ge 0$ as $R \to \infty$.  Define $\chi_i^R:=sgn(\zeta_i^R) \zeta_i^R$ and multiply system (\ref{zetar}) with $sgn(\zeta^R_i)$ to get that 
 \begin{equation}\label{chir}
  \left\{
                      \begin{array}{ll}
                                             \Delta \chi_i^R = H_{u_iu_i} \chi_i^R - \sum_{j\neq i} sgn (H_{u_iu_j}) H_{u_iu_j} \chi^R_j- \lambda_1(R) \chi_i^R , & \hbox{if $|x|< R$,} \\
                       \chi_i^R=0, & \hbox{if $|x|= R$.}
                                                                       \end{array}
                    \right.\end{equation}
Note that to get (\ref{chir}) we have used  (\ref{zetamono}), i.e., $sgn(\zeta_i^R)=-sgn(H_{u_iu_j})sgn(\zeta_j^R)$. Since now $\chi_i^R$ is a nonnegative solution for (\ref{chir}), Harnack's inequality yields that for any compact subset $K$, $\max_K |\chi_i^R|	\le C(K) \min_K |\chi_i^R| $ for all $i=1,\cdots,m$ with the latter constant being independent of $\chi_i^R$. Standard elliptic estimates also yield that the family $(\chi_i^R)_R$ have also
uniformly bounded derivatives on compact sets. It follows that for a subsequence $(R_k)_k$ going to infinity, $(\chi_i^{R_k})_k$  converges in $C^2_{loc} (\mathbf{R}^N)$ to some $\chi_i \in C^2(\mathbf {R}^N)$ and that $\chi_i\ge0$. From (\ref{chir}) we see that $\chi_i$ satisfies 
\begin{eqnarray}\label{chi}
 \Delta \chi_i&=& H_{u_iu_i} \chi_i - \sum_{j\neq i} sgn (H_{u_iu_j}) H_{u_iu_j} \chi_j- \lambda_1 \chi_i \\&\le&\nonumber  (H_{u_iu_i} - \lambda_1 )\chi_i 
  \end{eqnarray} 
Since $\chi_i\ge 0$ and $H_{u_iu_j}$ is bounded, the strong maximum principle yields that  either $\chi_i=0$ or $\chi_i>0$ in $\mathbf{R}^N$.  If now $\chi_i=0$, then  from (\ref{chi}) we have $ \sum_{j\neq i} sgn (H_{u_iu_j}) H_{u_iu_j} \chi_j=0$ which means $\chi_j=0$ if $j\neq i$, which contradicts  (\ref{normal}). It follows that $\chi_i>0$ for all $i=1,\cdots,m$. Set now $\phi_i:=sgn(\theta_i) \chi_i$ for $i=1,\cdots, m$ and observe that  $(\phi_i)_i$ satisfy (\ref{pointstability}) and that $H_{u_i,u_j} \phi_j \phi_i\le 0$ for $i<j$, which means that $u$ is a pointwise stable solution.  \hfill $ \Box$

\section{De Giorgi type results}

 We first establish a geometric Poincar\'{e} inequality for stable solutions of system (\ref{main}), which will enable us to get not only De Giorgi type results but also certain rigidity properties on the gradient of the solutions.

\begin{thm}\label{poin}
 Assume that  $m,N\ge 1$ and $\Omega\subset\mathbf{R}^N$ is an open set.  Then, for any $\eta=(\eta_k)_{k=1}^m \in C_c^1(\Omega)$, the following inequality holds for any classical stable solution $u \in C^2(\Omega)$ of (\ref{main})
\begin{eqnarray}\label{poincare}
\nonumber\sum_{i}   \int_{\Omega}|\nabla u_i|^2   |\nabla \eta_i|^2&\ge& \sum_{i} \int_{|\nabla u_i|\neq 0}\left(   |\nabla u_i|^2 \mathcal{A}_i^2 + | \nabla_T |\nabla u_i| |^2  \right)\eta_i^2\\&&+\sum_{i\neq j} \int_{\Omega}  \left( \nabla u_i \cdot\nabla  u_j \eta_i^2 - |\nabla u_i| |\nabla u_j| \eta_i \eta_j\right)H_{u_iu_j},
  \end{eqnarray} 
 where $\nabla_T$ stands for the tangential gradient along a given level set of $u_i$ and 
$\mathcal{A}_i^2$ for the sum of the squares of the principal curvatures of such a level set.
\end{thm}
\textbf{Proof:} Let $\eta=(\eta_1,...,\eta_m)$ and $\eta_i\in C^1_c(\Omega)$. Test the stability inequality (\ref{stability}) with $\zeta_i=|\nabla u_i|\eta_i$  to get 
 \begin{eqnarray}\label{stablepoin}
\nonumber 0&\le& \sum_{i} \int_{\Omega}  |\nabla (|\nabla u_i|\eta_i) ^2  +\sum_{i,j} \int_{\Omega} H_{u_iu_j} |\nabla u_i| |\nabla u_j|\eta_i\eta_j\\
&=&\nonumber   \sum_{i}   \int_{\Omega} |\nabla |\nabla u_i| |^2 {\eta_i^2}+\sum_{i}   \int_\Omega |\nabla \eta_i|^2|\nabla u_i|^2 +\frac{1}{2} \sum_{i}   \int_{\mathbf{R}^N} \nabla |\nabla u_i|^2\cdot\nabla {\eta_i^2} \\&&+ \sum_{i} \int_{\Omega}  H_{u_iu_i}|\nabla u_i|^2\eta_i^2+\sum_{i\neq j} \int_{\Omega} H_{u_iu_j} |\nabla u_i| |\nabla u_j|\eta_i\eta_j.
  \end{eqnarray} 
Differentiate the $i^{th}$ equation of (\ref{main}) with respect to $x_k$ for each $i=1,2,...,m$ and multiply with $\partial_k u$ to get 
$$ \partial_k u_i \Delta \partial_k u_i= \sum_{j} H_{u_iu_j} \partial_k u_j\partial_k u_i=H_{u_iu_j}|\partial_k u_i|^2+\sum_{j\neq i} H_{u_iu_j} \partial_k u_i \partial_k u_j.  $$
Multiply both sides with $\eta_i^2$ and integrate by parts to obtain
 \begin{eqnarray*}\label{}
\nonumber \int_{\Omega}  H_{u_iu_j}|\partial_k u_i|^2\eta_i^2&=&  \int_{\Omega}\partial_k u_i \Delta \partial_k u_i\eta_i^2- \sum_{j\neq i}  \int_{\mathbf{R}^N}  H_{u_iu_j} \partial_k u_i \partial_k u_j  \eta^2_i\\&=&-\int_{\Omega} |\nabla \partial_k u_i|^2\eta_i^2 - \frac{1}{2} \int_{\Omega} \nabla |\partial_k u_i|^2 \cdot\nabla \eta_i^2- \sum_{j\neq i}  \int_{\Omega}  H_{u_iu_j} \partial_k u_i \partial_k u_j  \eta_i^2.
  \end{eqnarray*} 
By summing over the index $k$,  we obtain
   \begin{eqnarray}\label{pde}
\nonumber  \int_{\Omega}  H_{u_iu_i}|\nabla u_i|^2\eta_i^2&=&- \sum_{k} \int_{\Omega} |\nabla \partial_k u_i|^2\eta_i^2 - \frac{1}{2} \int_{\Omega} \nabla |\nabla u_i|^2 \cdot\nabla \eta_i^2\\&&- \sum_{k} \sum_{j\neq i}  \int_{\Omega}  H_{u_iu_j} \partial_k u_i \partial_k u_j  \eta_i^2.
  \end{eqnarray} 
Combine (\ref{stablepoin})  and (\ref{pde}) to get 
\begin{eqnarray*}
\sum_{i}   \int_{\Omega}  |\nabla \eta_i|^2|\nabla u_i|^2 &\ge& \sum_{i} \int_{|\nabla u_i|\neq 0}\left(  \sum_{k} |\nabla \partial_k u_i|^2-| \nabla| \nabla u_i | |^2  \right)\eta_i^2\\&&+\sum_{i\neq j} \int_{\Omega} \left( \nabla u_i \cdot\nabla  u_j \eta_i^2 -|\nabla u_i| |\nabla u_j|\eta_i\eta_j \right)H_{u_iu_j}.
  \end{eqnarray*} 
 According to formula (2.1) given in  \cite{sz}, the following geometric identity between the tangential gradients and curvatures holds. For any $w \in C^2(\Omega)$
 \begin{eqnarray}\label{identity} \sum_{k=1}^{N} |\nabla \partial_k w|^2-|\nabla|\nabla w||^2=
\left\{
                      \begin{array}{ll}
                       |\nabla w|^2 (\sum_{l=1}^{N-1} \mathcal{\kappa}_l^2) +|\nabla_T|\nabla w||^2 & \hbox{for $x\in\{|\nabla w|>0\cap \Omega \}$,} \\
                       0 & \hbox{for $a.e.$ $x\in\{|\nabla w|=0\cap \Omega \}$,}
                                                                       \end{array}
                    \right.   \end{eqnarray} 
 where $ \mathcal{\kappa}_l$ are the principal curvatures of the level set of $w$ at $x$ and $\nabla_T$ denotes the orthogonal projection of the gradient along this level set . In light of this formula, we finally get (\ref{poincare}).
     
\hfill $ \Box$

 \begin{rem}
Note that for the case of $m=1$ the use of (\ref{identity}) and of $\zeta=|\nabla u|\eta $ in the stability (or semi-stability) condition (\ref{stability}) was first exploited by Sternberg and Zumbrun \cite{sz} to study semilinear phase transitions problems. Later on, Farina, Sciunzi, and Valdinoci \cite{fsv} used it to reprove the De  Giorgi's conjecture in dimension two, and Cabr\'{e}  used it (see Proposition 2.2 in \cite{c}) to prove the boundedness of extremal solutions of semilinear elliptic equations with Dirichlet boundary conditions on a convex domain up to dimension four.
\end{rem}

Here is an application of the above geometric Poincar\'{e} inequality for stable solutions of (\ref{main}).
 
\begin{thm} \label{superDG} Any bounded stable solution $u$ of an orientable system (\ref{main}) in $\mathbf{R}^2$ is one-dimensional. 
Moreover, if $H_{u_iu_j}$ is not identically zero, then for $i\neq j$,
\begin{equation}
\hbox{ $\nabla u_i=C_{i,j} \nabla u_j$ for all $x\in \mathbf{R}^2$,}
\end{equation}
 where   $C_{i,j}$ are constants whose sign is opposite to the one of $H_{u_iu_j}$.

\end{thm}
 
\noindent \textbf{Proof:} Fix the following standard test function 
$$\chi (x):=\left\{
                      \begin{array}{ll}
                        \frac{1}{2}, & \hbox{if $|x|\le\sqrt{R}$,} \\
                      \frac{ \log \frac{R}{|x|}}{{\log R}}, & \hbox{if $\sqrt{R}< |x|< R$,} \\
                       0, & \hbox{if $|x|\ge R$.}
                                                                       \end{array}
                    \right.$$
Since the system (\ref{main}) is {\it orientable}, there exist nonzero functions $\theta_k\in C^1(\mathbf{R}^N)$, $k=1,\cdots,m$, which do not change sign such that 
 \begin{equation}\label{sign}
  H_{u_iu_j}\theta_i\theta_j\le 0, \ \ \ \text{for all} \ \ i,j\in\{1,\cdots,m\} \ \text{and} \ \ i<j.
  \end{equation} 
Consider  $\eta_k:=sgn(\theta_k) \chi$ for $1\le k\le m$, where again $sgn (x)$ is the Sign function. The geometric  Poincar\'{e} inequality (\ref{poincare}) yields
\begin{eqnarray}\label{poincare2}
\nonumber  \int_{B_R\setminus B_{\sqrt {R} } } \sum_{i}  |\nabla u_i|^2   |\nabla \chi|^2&\ge& \sum_{i} \int_{|\nabla u_i|\neq 0}\left(   |\nabla u_i|^2 \mathcal{\kappa}_i^2 + | \nabla_T |\nabla u_i| |^2  \right)\chi^2\\&&\nonumber+\sum_{i\neq j} \int_{\mathbf{R}^N}  \left( \nabla u_i \cdot\nabla  u_j  - sgn(\theta_i) sgn(\theta_j) |\nabla u_i| |\nabla u_j|\right)H_{u_iu_j} \chi^2\\&=&I_1+I_2 .
  \end{eqnarray} 
Note that  $I_1$ is clearly nonnegative. Moreover,  
(\ref{sign}) yields that $H_{u_iu_j}sgn(\theta_i) sgn(\theta_j) \le 0$ for all $i<j$, and 
therefore, $I_2$ can be written as 
  $$I_2=\sum_{i\neq j} \int_{\mathbf{R}^N}  \left(sgn(H_{u_iu_j})  \nabla u_i \cdot\nabla  u_j  + |\nabla u_i| |\nabla u_j|\right)H_{u_iu_j} sgn(H_{u_iu_j}) \chi^2,
  $$
which is also nonnegative. 

On the other hand, since  $$\int_{B_R\setminus B_{\sqrt {R} } } \sum_{i}  |\nabla u_i|^2   |\nabla \chi|^2 \le C \left\{
                      \begin{array}{ll}
                        \frac{1}{\log R}, & \hbox{if $N=2$,} \\
                       \frac{R^{N-2}+ R^{(N-2)/2}}{|N-2||\log R|^2}, & \hbox{if $N\neq 2$,}
                                                                       \end{array}
                    \right.$$ one can see that in dimension two the left hand side of (\ref{poincare2}) goes to zero as $R\to \infty$. Since $I_1=0$, one concludes that all $u_i$ for $i=1,\cdots,m$ are one-dimensional and from the fact that $I_2=0$, provided $H_{u_iu_j}$ is not identically zero, we obtain that for all $x\in\mathbf{R}^2$, 
                    $$-sgn(H_{u_iu_j})  \nabla u_i \cdot\nabla  u_j  = |\nabla u_i| |\nabla u_j|,$$ 
                   which completes the proof of the theorem.   \hfill $ \Box$

    Now, we are ready to state and prove the main result of this paper.

\begin{thm} \label{Degiorgithm}
 Conjecture (\ref{conj2}) holds for $N\leq 3$.
\end{thm}

\textbf{Proof:} Let again $\phi_i := \partial _N u_i$ and $\psi_i:=\nabla u_i\cdot\eta$ for any fixed $\eta=(\eta',0)\in \mathbf{R}^{N-1}\times\{0\}$ in such a way 
that $\sigma_i:=\frac{\psi_i}{\phi_i}$  is a solution of system (\ref{div}) for $h_{i,j}(x)=H_{u_iu_j}\phi_i(x)\phi_j(x)$ and $f$ to be the identity. Since $|\nabla u_i|\in L^\infty(\mathbf{R}^N)$, we have $||\phi_i\sigma_i||_{L^\infty(\mathbf{R}^N)}< \infty$.

In dimension $N=2$,  assumption (\ref{liouassum}) holds and Proposition \ref{liouville} then yields that $\sigma_i$  is constant, which finishes the proof as argued before. 

 In dimension $N=3$, we shall follow ideas  used by Ambrosio-Cabr\'{e}  \cite{ac} and  Alberti-Ambrosio-Cabr\'{e} \cite{aac} in the case of a single equation. We first note that $u$ being $H$-monotone means that  $u$ is a stable solution of (\ref{main}). Moreover, the function $v(x_1, x_2):=\lim_{x_3\to \infty}u(x_1, x_2, x_3)$ is also a bounded stable solution for (\ref{main}) in $\mathbf{R}^2$. Indeed, it suffices to  
 test (\ref{stability}) on $\zeta_k(x)=\eta_k(x')\chi_{R} (x_N)$ where $\eta_k\in C_c^1(\mathbf{R}^{N-1})$ and  $\chi_R\in C_c^1(\mathbf{R})$ is defined as
 $$\chi_R (t):=\left\{
                      \begin{array}{ll}
                       1, & \hbox{if $R+1< t <2R+1$,} \\
                     0, & \hbox{if $t<R \ \ \text{or} \ \  t >2R+2$,} 
                                                                                             \end{array}
                    \right.$$
 for $R>1$, $0\le \chi_R\le 1$ and $0\le \chi'_R\le 2$. Note also that since $u$ is an $H$-monotone solution, the system (\ref{main}) is then orientable. It follows from Theorem \ref{superDG} that $v$ is one dimensional and consequently the energy of $v$ in a two-dimensional ball of radius $R$ is bounded by a multiple of $R$, which yields that
 \begin{equation}
\label{boundEt}
 \limsup_{t\to\infty} E(u^t) \le C R^{2},
 \end{equation}
 where here $u^t(x'):=u(x',x_n+t)$ for $t\in\mathbf{R}$ and $E_R(u)=\int_{B_{R}} \frac{1}{2}|\nabla u|^2 + H(u)-c_u d\mathbf x$ for $c_u:=\inf H(u)$.

To finish the proof, we shall show that  
\begin{equation}
\label{bound}
 \int_{B_{R}} |\nabla u|^2 \le C R^{2}.
 \end{equation}
Note that shifted function $u^t$ is also a bounded solution of (\ref{main}) with $|\nabla u_i^t| \in L^{\infty}(\mathbf{R}^N)$, i.e.,
\begin{eqnarray}
\label{maint}
 \Delta u^t&=& \nabla H(u^t)   \ \ \text{in}\ \ \mathbf{R}^N,
  \end{eqnarray}
and also 
  \begin{eqnarray}
\label{monot}
\partial_t u^t_i>0 > \partial_t u^t_j \ \
\ \ \text{for all $i\in I$ and $j\in J$ and in} \  \mathbf{R}^N.
  \end{eqnarray}
Since  $u_i^t$ converges to $v_i$ in $C^1_{loc}(\mathbf{R}^N)$ for all $i=1,\cdots,m$, we have $$ \lim_{t\to\infty} E(u^t)= E(v).
$$
Now, we claim that the  following upper bound for the energy holds. 
\begin{equation}\label{energyt}
E_R(u) \le  E_R(u^t) +M \int_{\partial B_{R}} \left(  \sum_{i\in I} (u_i^t-u_i)+ \sum_{j\in J}(u_j-u_j^t)  \right)  dS \ \ \text{for all}\ \ t \in\mathbf{R^+},
\end{equation}  
where $M=\max_i ||\nabla u_i||_{L^\infty(\mathbf{R}^N)}$. Indeed, by differentiating the energy functional along the path $u^t$, one gets 
\begin{eqnarray}\label{der-energy}
\partial_t E_R(u^t)= \int_{B_{R}} \nabla u^t\cdot \nabla (\partial _t u^t)+ \int_{B_{R}} \nabla H(u^t) \partial_t u^t,
    \end{eqnarray}
    where $\nabla H(u^t) \partial_t u^t=\sum_{i} H_{u_i}(u^t) \partial_t u^t_i$.
Now, multiply (\ref{maint}) with $\partial_t u^t$, to obtain 
  \begin{eqnarray}
\label{partialt}
\hfill - \int_{B_{R}} \nabla u^t\cdot \nabla (\partial _t u^t)+ \int_{\partial B_{R}} \partial_\nu u^t \partial_t u^t &=&\int_{B_{R}} \nabla  H(u^t) \partial_t u^t.
  \end{eqnarray}
 From (\ref{partialt}) and  (\ref{der-energy}) we obtain 
  \begin{eqnarray}\label{der-energy-bd}
\partial_t E_R(u^t)= \int_{\partial B_{R}} \partial_\nu u^t \partial_t u^t=\sum_{i}  \int_{\partial B_{R}} \partial_\nu u_i^t \partial_t u_i^t.
   \end{eqnarray}
 Note that $-M\le \partial_\nu u^t\le M$ and $\partial_t u_i^t>0>\partial_t u_j^t$ for $i\in I$ and $j\in J$ . Therefore, 
  \begin{eqnarray}\label{der-energy-bdM}
\partial_t E_R(u^t) \ge M \int_{\partial B_{R}}\left( \sum_{j} \partial_t u_j^t - \sum_{i}\partial_t u_i^t \right) d S.
   \end{eqnarray}
   On the other hand,
 \begin{eqnarray}\label{der-energy-last}
 E_R(u)&=& E_R(u^t)- \int_{0}^{t} \partial_t E_R(u^s) ds,\nonumber\\
 &\le&E_R(u^t) +M \int_{0}^{t}  \int_{\partial B_{R}} \left(  \sum_{i} \partial_s u_i^s-  \sum_{j} \partial_s u_j^s \right) dS ds\nonumber\\
 &=&E_R(u^t) +M \int_{\partial B_{R}}  \left(  \sum_{i} (u_i^t-u_i)+ \sum_{j}(u_j-u_j^t)  \right) dS.
 \end{eqnarray}
To finish the proof of the theorem just note that $u_i<u_i^t$ and $u_j^t<u_j$ for all $i\in I$, $j\in J$ and $t\in \mathbf{R^+}$. Moreover, from (\ref{boundEt}) we have $\lim_{t\to\infty} E_R(u^t)\le CR^2$. Therefore, (\ref{der-energy-last})  yields  $$ E_R(u) \le C |\partial B_R|\le C R^{2}, $$
and we are done.  \hfill  $\Box$\\

The above proof suggests that --just as in the case of a single equation-- any $H$-monotone solution $u$ of $(\ref{main})$   must satisfy the following estimate 
\begin{equation}\label{estimate}
\hbox{$ \int_{B_{R}} |\nabla u|^2 \le C R^{N-1}$ for any $R>1$,}
 \end{equation}
 for some constant $C>0$. This can be done in the following particular case.

\begin{thm} \label{energythm}
If $u$ is a bounded $H$-monotone solution of (\ref{main}) such that for 
 $i=1,..,m$,
$$\lim_{x_N\to \infty} u_i(\mathbf x',x_N)=a_i, \ \ \ \forall \mathbf x=(\mathbf x', x_N)\in\mathbf{R}^{N-1}\times\mathbf{R}$$
where $a_i$ are constants, then 
\begin{equation}\label{energybound}
E_R(u)=\int_{B_{R}} \frac{1}{2}|\nabla u|^2 + H(u)-H(a) d\mathbf x \le C R^{N-1},
\end{equation}
where  $a=\{a_i\}_{i=1}^{i=m}$ and $C$ is a positive constant  independent of $R$.
 \end{thm} 
\textbf{Proof:} We first  note the following decay on the energy of the shifted function $u^t$ as defined above, 
\begin{equation}\label{energydecay}
\lim_{t\to\infty} E_R(u^t)=0.
\end{equation}
Indeed, since $u^t$ is convergent to $a$ pointwise, one can see that $$\lim_{t\to\infty} \int_{B_{R}} (H(u^t) -H(a)) d\mathbf x \to 0.$$
Therefore, we need to prove that
$$\lim_{t\to\infty}  \int_{B_{R}}  |\nabla u_i^t|^2 d\mathbf x \to 0.$$
To do so, multiply both sides  of (\ref{maint}) with $u_i^t-a_i$ and integrate by parts to get
\begin{eqnarray*}
\label{}
\hfill -\int_{B_{R}}  |\nabla u_i^t|^2+ \int_{\partial B_{R}} \partial_\nu u_i^t (u_i^t-a_i) &=&\int_{B_{R}}  \nabla H(u^t) (u_i^t-a_i),
  \end{eqnarray*}
which yields  (\ref{energydecay}).

To get the energy bound in (\ref{energybound}), one can follow the proof of the previous theorem to end up with 
 \begin{eqnarray*}
  E_R(u) \le E_R(u^t) +C  |\partial B_R| \ \ \ \text{for all} \ \ t\in \mathbf{R^+}.
 \end{eqnarray*}
To conclude, it suffices to send $t\to\infty$ and to use the fact that $\lim_{t\to\infty} E_R(u^t)=0$ to finally obtain that  $$ E_R(u) \le C |\partial B_R|\le C R^{N-1}. $$ 

\begin{rem} \rm
 Using Pohozaev type arguments one can see that 
\begin{equation}
\hbox{$\Gamma_R=\frac{E_R(u)}{R^{N-1}}$ is increasing}
\end{equation}
 provided the following pointwise estimate holds:
\begin{equation}
|\nabla u|^2 \le 2 H(u).
\end{equation}
Note that this is an extension of the pointwise estimate that Modica \cite{m} proved in the case of a single equation. It  is  still not known for systems, though Caffarelli-Lin in \cite{cl} and later, Alikakos in \cite{ali} have shown, in the case where  $H\ge 0$, the following 
weaker monotonicity formula, namely that 
\begin{equation}
\hbox{$\Lambda_R=\frac{E_R(u)}{R^{N-2}}$ is increasing in $R$.}
\end{equation}
\end{rem} 

\begin{rem} \rm The $H$-monotonicity assumption seems to be crucial for concluding that the solutions are one-dimensional. Indeed, it was shown in \cite{ABG} that when $H$ is a multiple-well potential on  $\mathbf{R}^2$, the system has entire ÒheteroclinicÓ solutions $(u, v)$, meaning that for each fixed $x_2 \in \mathbf{R}$, they connect (when $x_1 \to \pm \infty$) a pair of constant global minima of W, while if $x_2 \to \pm \infty$, they connect
a pair of distinct one dimensional stationary wave solutions $z_1(x_1)$ and $z_2(x_1)$. Note that these convergence are even uniform, which means that the corresponding Gibbons conjecture for systems of equations is not valid in general,  without the assumption of $H$-monotonicity. 

\end{rem}

 \end{document}